\documentclass[twoside,leqno,10pt]{amsart}
\usepackage{amsfonts}
\usepackage{amsmath}
\usepackage{amscd}
\usepackage{amssymb}
\usepackage{amssymb}
\usepackage{amsthm}
\usepackage{amsrefs}
\everymath={\displaystyle}
\begin{document}
\newtheorem{theorem}{Theorem}
\newtheorem{lemma}{Lemma}
\newtheorem{corollary}{Corollary}
\newtheorem{conjecture}{Conjecture}
\numberwithin{equation}{section}
\newcommand{\dif}{\mathrm{d}}
\newcommand{\intz}{\mathbb{Z}}
\newcommand{\ratq}{\mathbb{Q}}
\newcommand{\natn}{\mathbb{N}}
\newcommand{\comc}{\mathbb{C}}
\newcommand{\rear}{\mathbb{R}}
\newcommand{\prip}{\mathbb{P}}
\newcommand{\uph}{\mathbb{H}}

\title{\bf Large Sieve Inequality with Quadratic Amplitudes}
\author{Liangyi Zhao}
\maketitle

\begin{abstract}
In this paper, we develop a large sieve type inequality with quadratic amplitude.  We use the double large sieve to establish non-trivial bounds.\end{abstract}

\noindent{\bf keywords:} large sieve, quadratic amplitudes \newline
{\bf 2000 Mathematics Classification:} 11B57, 11L07, 11L40, 11N35

\section{Introduction and Statements of the Results}

The large sieve was an idea originated by Yu. V. Linnik \cite{JVL1} in 1941.  He also made application to distributions of quadratic non-residues.  Since then, the idea has been refined and perfected by many. \newline 

We denote $\| x \| = \min_{k \in \intz} |x-k|$ for $x \in \rear$.  A set of real numbers $\{ x_k \}$ is said to be $\delta$-spaced modulo 1 if $\| x_j-x_k \| > \delta$, for any $j \neq k$. \newline

The large sieve inequality, which we henceforth refer to as the classical large sieve inequality, is stated as follows.  Different elegant proofs of the theorem can be found in \cites{HD, EB1, EB2, PXG, HM2, HM, HLMRCV}.  The theorem, in the following form, was first introduced by Davenport and Halberstam, \cite{DH1}.

\begin{theorem} \label{classicls}
Let $\{ a_n \}$ be an arbitrary sequence of complex numbers, $\{ x_k \}$ be real numbers that are $\delta$-spaced modulo 1, and $M$, $N$ be integers with $N>0$, then we have
\begin{equation} \label{classiclseq}
\sum_k \left| \sum_{n=M+1}^{M+N} a_n e \left( x_k n \right) \right|^2 \ll \left( \delta^{-1} +N \right) Z,
\end{equation}
where the implied constant is absolute and henceforth we set
\[ Z = \sum_{n=M+1}^{M+N} |a_n|^2 . \]
\end{theorem}
Save for the computation of the implied constant, the above inequality is the best possible.  Montgomery and Vaughan \cite{HLMRCV} proved that the ``$\ll$'' in \eqref{classiclseq} can be replaced by ``$\leq$.''  Moreover, Paul Cohen and A. Selberg have shown independently that 
\begin{equation*}
\sum_k \left| \sum_{n=M+1}^{M+N} a_n e \left( x_k n \right) \right|^2 \leq \left( \delta^{-1} -1+N \right) Z.
\end{equation*}
which is absolutely the best possible, since Bombieri and Davenport \cite{BD1} gave examples of $\{x_k\}$ and $a_n$, with $\delta \to 0$, $N \to \infty$ such that $N\delta \to \infty$ and
\begin{equation*}
\sum_k \left| \sum_{n=M+1}^{M+N} a_n e \left( x_k n \right) \right|^2 = \left( \delta^{-1} -1+N \right) Z.
\end{equation*}

Theorem~\ref{classicls} admits corollaries for additive and multiplicative characters.  In short, we have
\begin{equation} \label{lsadd}
\sum_{q=1}^Q \sum_{\substack{a \; \bmod{\; q} \\ \gcd(a,q)=1}} \left| \sum_{n=M+1}^{M+N} a_n e \left( \frac{a}{q} n \right) \right|^2 \ll (Q^2+N) Z.
\end{equation}

Moreover, for multiplicative characters, we have
\begin{equation} \label{lsmult}
\sum_{q=1}^Q \frac{q}{\varphi(q)} \sideset{}{^{\star}}\sum_{\chi
\; \bmod{ \;q}} \left| \sum_{n=M+1}^{M+N} a_n \chi (n) \right|^2 \ll (Q^2+N) Z
\end{equation}
where here and after, $\sideset{}{^{\star}}\sum$ means that the sum runs over primitive characters modulo a specified modulus only.  Also, as usual, $\varphi(q)$ is the number of primitive residue classes modulo $q$, {\it id est} the Euler $\varphi$ function.  The cases in which the outer summation in \eqref{lsadd} and \eqref{lsmult} run over only a set of special characters are investigated in \cites{LZ1, LZ4}. \newline

In this paper, we shall be interested in obtaining bounds of following form
\begin{equation*}
 \sum_{k=1}^K \left| \sum_{n=M+1}^{M+N} a_n e \left( x_k f(n) \right) \right|^2 \ll \Delta(\delta,N,\epsilon) Z,
\end{equation*}
where $\{ x_k \}$ are some well-spaced real numbers, and $f(x) = \alpha x^2 + \beta x + \gamma$ with $\alpha$, $\beta$, $\gamma \in \mathbb{R}$ and $\alpha >0$. \newline

We prove the following.

\begin{theorem} \label{lssfreq5}
Let $F(Q) = \{ p/q \in \ratq : 0 \leq p < q, \; \gcd(p,q)=1 \}$, and $f(x)=\alpha x^2 + \beta x + \gamma$ with $\beta/\alpha = a/b \in \ratq$, with $b>0$ and $\gcd(a,b)=1$.  We have
\begin{equation} 
\sum_{x \in F(Q)} \left| \sum_{n=M+1}^{M+N} a_n e \left( x f(n) \right) \right|^2 \ll (Q^2+Q\sqrt{\alpha N(|M|+N+a/b)+1}) \Pi Z,
\end{equation}
where the implied constant depends on $\epsilon$ alone and
\[ \Pi = \left( \frac{b}{\alpha}+1 \right)^{\frac{1}{2}+\epsilon} [Nb(|M|+N)+|a| + b/\alpha ]^{\epsilon}. \]
\end{theorem}

The above theorem may be extended to the case $\frac{\beta}{\alpha} \in \mathbb{R}$ by approximating $\frac{\beta}{\alpha}$ with a rational number and the proof works almost exactly the same.  More precisely, we have

\begin{theorem} \label{prop}
The result of Theorem~\ref{lssfreq5} still holds if the condition $\frac{\beta}{\alpha} \in \ratq$ is replaced by $\frac{\beta}{\alpha} \in \mathbb{R}$ and there exists $\frac{a}{b} \in \ratq$ with $\gcd(a,b)=1$ and
\begin{equation*}
\left| \frac{\beta}{\alpha} - \frac{a}{b} \right| < \frac{1}{4bN}.
\end{equation*}
\end{theorem}

\section{Preliminary lemmas}

We begin by quote the duality principle which says the norm of a bounded linear operator in a Banach space is the same as that of its adjoint operator.  To us, it amounts to the changing of orders of summations.

\begin{lemma}[Duality Principle]\label{dual}
Let $T=[t_{mn}]$ be a square matrix with entries from the complex numbers.  The following two statements are equivalent:
\begin{enumerate}
\item For any absolutely square summable sequence of complex numbers $\{ a_n \}$, we have
\begin{equation*}
\sum_m \left| \sum_n a_n t_{mn} \right|^2 \leq D \sum_n |a_n|^2.
\end{equation*}
\item For any absolutely square summable sequence of complex numbers $\{ b_n \}$, we have
\begin{equation*}
\sum_n \left| \sum_m b_m t_{mn} \right|^2 \leq D \sum_m |b_m|^2.
\end{equation*}
\end{enumerate}
\end{lemma}

\begin{proof} See for example Theorem 288 of \cite{GHHJELGP}. \end{proof}

For completeness, as we shall mention the following lemma later on, we give

\begin{lemma} [Sobolev-Gallagher] \label{sg}
Let $f(x) : [a,b] \rightarrow \mathbb{C}$ be a smooth function with continuous first derivative on $(a,b)$.  Then given $u \in [a,b]$, we have
\[ |f(u)| \leq \frac{1}{b-a} \int_a^b |f(x)| \dif x + \int_a^b |f'(x)| \dif x, \; \mbox{and} \; \left| f \left( \frac{a+b}{2} \right) \right| \leq \frac{1}{b-a} \int_a^b |f(x)| \dif x + \frac{1}{2} \int_a^b |f'(x)| \dif x. \]
\end{lemma}

\begin{proof} This is quoted from \cite{HM} and of course from \cite{PXG}. \end{proof}

We shall make extensive use of the following lemma.

\begin{lemma}[Double Large Sieve] \label{dlsieve}
Suppose $x_1$, $\cdots$, $x_M$ and $y_1$, $\cdots$, $y_N$ are real numbers with
\[ -\frac{X}{2} \leq x_m \leq \frac{X}{2}, \; -\frac{Y}{2} \leq y_n \leq \frac{Y}{2}, \]
for $m=1$, $\cdots$, $M$, $n=1$, $\cdots$, $N$.  Put $\epsilon = X^{-1}$, $\delta = \frac{X}{XY+1}$, and $\Lambda (x) = \max (1-|x|,0)$.
Then we have
\begin{equation*}
\left| \sum_{m=1}^M \sum_{n=1}^N a_m b_n e \left( x_m y_n\right) \right|^2 \leq \left( \frac{\pi}{2} \right)^{4} A(\delta) B(\epsilon) (X Y +1),
\end{equation*}
where
\[ A(\delta) = \sum_{m=1}^M \sum_{r=1}^M |a_m||a_r| \Lambda \left( \frac{x_m-x_r}{\delta} \right) \; \mbox{and} \; B(\epsilon) = \sum_{n=1}^N \sum_{r=1}^N b_n \overline{b}_r \Lambda \left( \frac{y_n-y_r}{\epsilon} \right). \]
\end{lemma}
\begin{proof} This is the one-dimensional version of Lemma 5.6.6 from \cite{MNH} and has its origin in Lemma 2.4 in \cite{BI1}. \end{proof}

\section{Heuristics and Trivial Bounds}

We are interested in estimating the size of the following expression.  
\begin{equation*}
\sum_{k=1}^K
\left| \sum_{n=M+1}^{M+N} a_n e \left( x_k f(x) \right) \right|^2.
\end{equation*}

By the virtue of the duality principle, Lemma~\ref{dual}, it is equivalent to estimate for the following sum.
\begin{equation*}
\sum_{n=M+1}^{M+N} \left| \sum_{k=1}^K c_k e \left( x_k f(n) \right) \right|^2,
\end{equation*}
where $\{ c_k \}$ is an arbitrary sequence of complex numbers.  Therefore, if $f(n)=n^2$, then the classical large sieve inequality, in its dual form, gives
\begin{equation} \label{lsftriv}
\sum_{n=M+1}^{M+N} \left| \sum_{k=1}^K c_k e \left( x_k n^2 \right) \right|^2 \ll (\delta^{-1} + (|M|+N)^2 ) \sum_{k=1}^K |c_k|^2,
\end{equation}
where the $N^2$ in the \eqref{lsftriv} comes from changing the variables and then filling in all the non-squares $n$'s.  Consequently, we have the corresponding result.
\begin{equation*}
\sum_{k=1}^K \left| \sum_{n=M+1}^{M+N} a_n e(x_k n^2) \right|^2 \ll (\delta^{-1}+(|M|+N)^2) Z.
\end{equation*}

In the light of how we arrived at \eqref{lsftriv}, one may tend to think that we should have the same estimate that the regular large sieve inequality has.  Namely, the upper bound should be, essentially,
\begin{equation} \label{goal}
(\delta^{-1} +N)Z.
\end{equation}

However, the bound in \eqref{goal} is false in general.  Let $\{ x_k \}$ be the Farey sequence of order $Q$.  Therefore, $\delta^{-1}=Q^2$ in this case.  Moreover, take $Q=p^2$ be a prime square and $a_n$ be $p$ if $n$ is a multiple of $p$ and zero otherwise.  Also take $M=0$ and $N$ a positive multiple of $p$.  Therefore, we have
\begin{equation} \label{countex1}
\sum_{q=1}^Q \sum_{\substack{a \; \bmod{ \;q} \\ \gcd(a,q)=1}} \left| \sum_{n=M+1}^{M+N} a_n e \left( \frac{a}{q} n^2 \right) \right|^2  \gg \sum_{\substack{a \; \bmod{ \;Q=p^2} \\ \gcd(a,Q)=1}} \left| p \frac{N}{p} \right|^2 \gg QN^2.
\end{equation}
where we minorize the sum over $q$ by the single term corresponding to $q=Q=p^2$.  But on the other hand,
\begin{equation} \label{countex2}
(\delta^{-1}+N) Z \asymp (Q^2+N) \frac{N}{p}
p^2 = Q^{\frac{5}{2}}N+Q^{\frac{1}{2}}N^2.
\end{equation}

Now the lower bound in \eqref{countex1} exceeds \eqref{countex2} whenever $N \gg Q^{3/2}$.  Thus we have a contradiction, and the result in \eqref{goal} is not attainable in general. \newline

Furthermore, if one consider the general case of $f(n) = \alpha n^2 + \beta n + \gamma$, a quadratic polynomial with real coefficients, one could mimic Gallagher's proof of the regular large sieve inequality in \cite{PXG}.  More precisely, we apply the Sobolev-Gallagher's Lemma, Lemma~\ref{sg}, to 
\[ \left( \sum_{n=M+1}^{M+N} a_n e \left( x f(n) \right) \right)^2 \]
on the intervals $[ x_k-\delta, x_k+\delta ]$ and then sum over the $x_k$'s.  We then get an upper bound for the sum of our interest,
\[ \delta^{-1} \int_0^1 \left| \sum_{n=M+1}^{M+N} a_n e (x f(n)) \right|^2 \dif x + 2 \pi \int_0^1 \left| \sum_{n=M+1}^{M+N} a_n e (x f(n)) \right| \left| \sum_{n=M+1}^{M+N} a_n f(n)e (x f(n)) \right| \dif x. \]

Applying Parseval's identity to the first integral, we infer that it is $\ll Z$.  By the virtues of H\"{o}lder's inequality to the second integral and then Parseval's inequality, we see that the second integral is majorized by
\[ \leq \left( \int_0^1 \left| \sum_{n=M+1}^{M+N} a_n e (x f(n))
\right|^2 \dif x \right)^{\frac{1}{2}} \left( \sum_{n=M+1}^{M+N} \left| a_n f(n) \right|^2 \right)^{\frac{1}{2}}. \]
Applying Parseval's identity again to the first factor, we get that above is bounded by
\begin{equation*}
\ll [\delta^{-1} + |\alpha| (|M|+N)^2] Z.
\end{equation*}
Theorem~\ref{lssfreq5} is better than the above if $Q \ll N$.

\section{Proof of Results}

We need the following lemma concerning the difference of the difference of quadratic polynomials.

\begin{lemma} \label{diffofdiff}
Let $S = [M+1, M+N] \cap \mathbb{Z}$, $\alpha>0$, and $\epsilon >0$ be given.  Set 
\[ g(x,y)=(x-y)(x+y+\frac{a}{b}), \]
with $a$, $b \in \mathbb{Z}$, $b>0$ and $\gcd(a,b)=1$.  For fixed $m$, $n \in S$, let $T$ denote the number of pairs of $m'$, $n' \in S$ with $g(m',n') \neq 0$ and
\begin{equation} \label{diffofdiff2}
| g(m,n)-g(m',n') | \leq \frac{1}{2\alpha}.
\end{equation}
Then we have
\begin{equation}
T \ll \left( \frac{b}{\alpha}+1 \right) \left[ Nb(|M|+N)+|a| + b/\alpha \right]^{\epsilon},
\end{equation}
where the implied constant depends only on $\epsilon$.
\end{lemma}

\begin{proof} First we note that \eqref{diffofdiff2} holds if and only if
\begin{equation} \label{diffofdiff4}
|bg(m,n)-bg(m',n')| =|(m-n)(bm+bn+a)-(m'-n')(bm'+bn'+a)| \leq \frac{b}{2\alpha}.
\end{equation}

For fixed $m$, $n \in S$, there are at most $b/\alpha +1$ integers that are no more than $b/(2\alpha)$ away from $bg(m,n)$.  If $k \neq 0$ is one such integer and $bg(m',n')=k$, then we can find $u$, $v \in \intz$ with $k=uv$\[ m'-n'=u \; \mbox{and} \; bm'-bn'+a=v. \]
We then see that 
\[ m'=\frac{bu+v-a}{2b} \; \mbox{and} \; n'=\frac{-bu+v+a}{2b}. \]
In other words, $m',n'$ are completely determined by the choice of $u$ and $v$.  Therefore, the number of pairs $(m',n')$ with $bg(m',n')=k\neq 0$ does not exceed 
\[ \tau(|k|) \ll |k|^\epsilon \ll (Nb(|M|+N+|a|)+b/\alpha)^\epsilon, \]
where $\tau(|k|)$ denotes the number of divisors of $|k|$. As we have already noted, the number of such $k$'s is at most $b/ \alpha+1$. Hence we have our desired result.
\end{proof}

Now we are ready to prove Theorem~\ref{lssfreq5}.

\begin{proof}(of Theorem~\ref{lssfreq5})
It is certainly elementary to write 
\[ f(x)= \alpha \left( x + \frac{\beta}{2\alpha} \right)^2 - \frac{\beta^2}{4\alpha} + \gamma \]
and hence it suffices to consider the sum
\[ \sum_{x\in F(Q)} \left| \sum_{s=M+1}^{M+N} c_s e \left( \alpha x \left( s + \frac{\beta}{2\alpha} \right)^2 \right) \right|^2, \]
with $\beta/\alpha = a/b$, $b>0$, $\gcd(a,b)=1$ and an arbitrary complex sequence $\{ c_s \}$.  Using the notations of Lemma~\ref{dlsieve}, we take
\[ x \in F(Q), \; y=g(s,t), \; X=2\alpha, \; Y=2|M|N+N^2+N\frac{a}{b}, \]
\[ a_m=1, \; \mbox{for all} \; m; \; b_n=c_s \overline{c}_t, \epsilon = \frac{1}{2\alpha}, \; \delta = \frac{2}{2(2|M|N+N^2+Na/b)+\alpha^{-1}}. \]
Recall that $g(s,t)$ is as defined in Lemma~\ref{diffofdiff}.  With these notations, we have
\begin{eqnarray*}
A(\delta) & = & \sum_x \sum_{x'} \Lambda \left[ \frac{4|M|N+2N^2+2Na/b+\alpha^{-1}}{2} \alpha (x-x') \right] \\
 & = & Q^2 + \sum_x \sum_{x' \neq x} \Lambda \left[ \frac{4|M|N+2N^2+2Na/b+\alpha^{-1}}{2} \alpha (x-x') \right].
\end{eqnarray*}
A summand in the above double sum is non-zero if and only if
\begin{equation} \label{difffarey}
 | x-x'| < \frac{2}{\alpha(4|M|N+2N^2+2Na/b)+1}
\end{equation}
Since $x$ and $x'$ are Farey fractions with denominator not exceeding $Q$, for each fixed $x$, the number of $x'$ satisfying \eqref{difffarey} is majorized as
\[ \ll \frac{Q^2}{\alpha N(|M|+N+a/b)+1} +1. \]
We therefore have
\[ A(\delta) \ll Q^2 + \frac{Q^4}{\alpha N(|M|+N+a/b)+1}. \]
It still remains to estimate $B(\epsilon)$.  We have
\begin{eqnarray}
\nonumber B(\epsilon) & = & \sum_{s,t} \sum_{s',t'} c_s \overline{c_t} \overline{c_{s'}} c_{t'} \Lambda(2 \alpha (g(s,t)-g(s',t')) \\
\label{Bestimate} & = &\sum_{s,t} \sum_{\substack{s',t' \\ g(s',t') =0}} c_s \overline{c_t} \overline{c_{s'}} c_{t'} \Lambda(2 \alpha g(s,t)) + \sum_{s,t} \sum_{\substack{s',t' \\ g(s',t') \neq 0}} c_s \overline{c_t} \overline{c_{s'}} c_{t'} \Lambda(2 \alpha (g(s,t)-g(s',t'))
\end{eqnarray}

Note that for a fixed $s$, the number of $t$'s with $g(s,t)=0$ is $O(1)$.  Hence the first term in \eqref{Bestimate} is 
\[ \ll \left( \sum_{s'} |c_{s'}|^2 \right) \left( \sum_{s,t} |c_s \overline{c_t}| \Lambda(2 \alpha g(s,t)) \right) \ll \left( \sum_{s'} |c_{s'}|^2 \right) \left( \sum_{s} |c_s|^2 + \sum_{\substack{s,t \\ 0 < |g(s,t)|<(2\alpha)^{-1}}} |c_s \overline{c_t}| \right), \]
by the virtue of the parallelogram rule.  As noted in Lemma~\ref{diffofdiff}, $0<|g(s,t)|<(2\alpha)^{-1}$ holds if and only if 
\begin{equation} \label{Bestimate2}
0 < |(s-t)(bs+bt+a)| < \frac{b}{2\alpha}.
\end{equation}
For a fixed $k \in \intz$ with $0<|k|<b/(2\alpha)$ and $s \in \intz$ fixed, the number of $t$'s with $t \in \intz$ and
\[ (s-t)(bs+bt+a) = k \]
is at most $\tau(|k|)$ by the same arguments in the proof of Lemma~\ref{diffofdiff}.  Hence with $s$ fixed, the number of $t$'s satisfying \eqref{Bestimate2} is, for any given $\eta >0$,
\[ \ll \left( \frac{b}{2\alpha} \right)^{\eta} \left( \frac{b}{\alpha} +1 \right) \ll \left( \frac{b}{\alpha} +1 \right)^{1+\eta}, \]
with the implied constant depending only on $\eta$.  Hence, the first term in \eqref{Bestimate} is estimated as
\[ \ll  \left( \sum_{s'} |c_{s'}|^2 \right) \left( \sum_{s}|c_s|^2 +  \left( \frac{b}{\alpha} +1 \right)^{1+\eta} \sum_s |c_s|^2 \right) \ll \left( \frac{b}{\alpha} +1 \right)^{1+\eta}  \left( \sum_{s'} |c_{s'}|^2 \right)^2. \]

We now use Lemma~\ref{diffofdiff} for the second term in \eqref{Bestimate}.  The summand of interest is zero unless
\begin{equation} \label{Bestimate3}
|g(s,t)-g(s',t')| < \frac{1}{2\alpha}.
\end{equation}
For fixed $s$ and $t$, the number of pairs $s'$ and $t'$ satisfying \eqref{Bestimate3} is
\[ \ll \left( \frac{b}{\alpha} +1  \right) \left( Nb \left( |M| +N \right) + |a| + \frac{b}{2\alpha} \right)^{\eta} . \]
This gives that the second term in \eqref{Bestimate} is 
\[ \ll \left( \frac{b}{\alpha} +1 \right) \left( Nb(|M|+N)+|a| + \frac{b}{2 \alpha} \right)^{\eta} \sum_{s,t} |c_s c_t|^2 = \left( \frac{b}{\alpha} +1 \right) \left( Nb(|M|+N)+|a| + \frac{b}{2 \alpha} \right)^{\eta} \left( \sum_{s} |c_s|^2 \right)^2 , \]
where the implied constant depends only on $\eta$.  Now combining everything and taking the square root, we have the desired result. \end{proof}

We prove Theorem~\ref{prop} next.

\begin{proof}(of the Theorem~\ref{prop}) The proof goes precisely the same as that of Theorem~\ref{lssfreq5}.  \eqref{Bestimate2} becomes
\[ 0 < |(s-t)(bs+bt+a)+(s-t)bz| < \frac{b}{2\alpha}, \]
where
\[ |z| < \frac{1}{4bN} \; \mbox{and hence} \; |(s-t)bz| < \frac{1}{4} . \]
Hence it suffices to estimate the number of pair $s$, $t$ such that
\[ 0 < |(s-t)(bs+bt+a)| < \frac{b}{2\alpha} + \frac{1}{4}. \]
Similarly, \eqref{diffofdiff4} becomes
\[ |(m-n)(bm+bn+a)-(m'-n')(bm'+bn'+a)+(m-n)bz-(m'-n')bz| \leq \frac{b}{2}, \]
and it suffices to estimate number of pairs $m'$ and $n'$ such that
\[ |(m-n)(bm+bn+a)-(m'-n')(bm'+bn'+a)| \leq \frac{b}{2} + \frac{1}{2}. \]
Now the arguments go precisely the same as those in the proofs of Theorem~\ref{lssfreq5} and Lemma~\ref{diffofdiff} and we arrive at the same results.
\end{proof}

\section{Notes}

With the counterexample given in Section 3, the author suspects the the result with $\Delta=Q^2+QN$ is essentially the best possible.  The author also believes that the method used in the chapter may be used to deal with large sieve inequality of higher power amplitudes.  Furthermore, we have been restricting our attention to only Farey fractions.  However, the result of double large sieve should enable us to consider any set of well-spaced real numbers. \newline

Finally, we note that the result of this paper in the case when $\beta/\alpha$ is irrational has been improved recently by S. Baier \cite{SB2}.  The situation with higher power amplitudes is also studied in the same paper.

\section*{Acknowledgments}

The author wishes to thank his thesis adviser, Professor Henryk Iwaniec, who suggested this problem to me and, of his support and advice, has been most generous.  This paper was finalized when the author was supported by a postdoctoral fellowship at the University of Toronto.

\bibliography{biblio}

% \bib, bibdiv, biblist are defined by the amsrefs package.
\begin{bibdiv}
\begin{biblist}

\bib{SB2}{article}{
    author={Baier, S.},
     title={The large sieve with quadratic amplitude},
      date={2005},
      note={pre-print},
}

\bib{EB1}{article}{
    author={Bombieri, E.},
     title={On the large sieves},
      date={1965},
   journal={Mathematika},
    volume={12},
     pages={201\ndash 225},
}

\bib{EB2}{article}{
    author={Bombieri, E.},
     title={Le grand crible dans la th\'eorie analytique des nombres},
      date={1974},
   journal={Soci\'et\'e Mathematics France},
    volume={18},
}

\bib{BD1}{article}{
    author={Bombieri, E.},
    author={Davenport, H.},
     title={Some inequalities involving trigonometrical polynomials},
      date={1969},
   journal={Annali Scuola Normale Superiore - Pisa},
    volume={23},
     pages={223\ndash 241},
}

\bib{BI1}{article}{
    author={Bombieri, E.},
    author={Iwaniec, H.},
     title={On the order of $\zeta \left( \frac{1}{2} +it \right)$},
      date={1986},
   journal={Annali Scuola Normale Superiore - Pisa},
    volume={XIII},
    number={3},
     pages={449\ndash 472},
}

\bib{HD}{book}{
    author={Davenport, H.},
     title={Multiplicative Number Theory},
   edition={Third Edition},
    series={Graduate Texts in Mathematics},
 publisher={Springer-Verlag},
   address={Barcelona, etc.},
      date={2000},
    volume={74},
}

\bib{DH1}{article}{
    author={Davenport, H.},
    author={Halberstam, H.},
     title={The values of a trigonometric polynomial at well spaced points},
      date={1966},
   journal={Mathematika},
    volume={13},
     pages={91\ndash 96},
      note={{\it Corrigendum and Addendum}, Mathematika {\bf 14} (1967),
  232-299},
}

\bib{PXG}{article}{
    author={Gallagher, P.~X.},
     title={The large sieve},
      date={1967},
   journal={Mathematika},
    volume={14},
     pages={14\ndash 20},
}

\bib{GHHJELGP}{book}{
    author={Hardy, G.~H.},
    author={Littlewood, J.~E.},
    author={P\'olya, G.},
     title={Inequalities},
 publisher={Cambridge University Press},
      date={1964},
}

\bib{MNH}{book}{
    author={Huxley, M.~N.},
     title={Area, Lattice Points and Exponential Sums},
    series={London Mathematical Society Monographs New Series},
 publisher={Clarendon Press},
   address={Oxford},
      date={1996},
    volume={13},
}

\bib{JVL1}{article}{
    author={Linnik, J.~V.},
     title={The large sieve},
      date={1941},
   journal={Doklady Akad. Nauk SSSR},
    volume={36},
     pages={119\ndash 120},
}

\bib{HM}{book}{
    author={Montgomery, H.~L.},
     title={Topics in Multiplicative Number Theory},
    series={Lecture Notes in Mathematics},
 publisher={Spring-Verlag},
   address={Barcelona, etc.},
      date={1971},
    volume={227},
}

\bib{HM2}{article}{
    author={Montgomery, H.~L.},
     title={The Analytic Priciples of Large Sieve},
      date={1978July},
   journal={Bulletin of the American Mathematical Society},
    volume={84},
    number={4},
     pages={547\ndash 567},
}

\bib{HLMRCV}{article}{
    author={Montgomery, H.~L.},
    author={Vaughan, R.~C.},
     title={The large sieve},
      date={1973},
   journal={Mathematika},
    volume={20},
     pages={119\ndash 134},
}

\bib{LZ1}{article}{
    author={Zhao, L.},
     title={Large sieve inequality for characters to square moduli},
      date={2004},
   journal={Acta Arithmetica},
    volume={112},
    number={3},
     pages={297\ndash 308},
}

\bib{LZ4}{article}{
    author={Zhao, L.},
     title={Large sieve inequality for special characters to prime square
  moduli},
      date={2004},
   journal={Functiones et Approximatio Commentarii Mathematici},
    volume={XXXII},
     pages={99\ndash 106},
}

\end{biblist}
\end{bibdiv}
\bibliographystyle{amsxport}

\noindent Department of Mathematics, University of Toronto \newline
40 Saint George Street, Toronto, ON M5S 2E4 Canada \newline
Email: {\tt lzhao@math.toronto.edu}

\end{document}